\documentclass[11pt]{amsart}
\usepackage{amsbsy}
\usepackage{amsmath}
\usepackage{amsbsy}
\usepackage{graphicx}
\usepackage{enumerate}
\usepackage[cp1250]{inputenc}
\usepackage{indentfirst}
\usepackage{latexsym}
\usepackage{amssymb}
\usepackage{amsfonts}
\usepackage{color}
\usepackage{dsfont}
\usepackage{mathtools}
\usepackage{amsthm}

\usepackage{tikz}
\usetikzlibrary{arrows}
\usepackage{xcolor}
\usepackage{xypic}
\usepackage{cite}
\usepackage{float}
\usepackage{filecontents}
\usepackage{graphicx}
\usepackage{epstopdf}
\usepackage{epsfig}

\definecolor{darkgreen}{rgb}{0,.7,.3}

\theoremstyle{definition}
\newtheorem{definition}{Definition}[section]
\newtheorem{remark}[definition]{Remark}

\theoremstyle{plain}
\newtheorem{lemma}[definition]{Lemma}

\newtheorem{proposition}[definition] {Proposition}

\newcommand{\Lab}{\text{\rm \bf Lab}}

\newcommand{\Cat}{\text{\rm  Cat}}
\newcommand{\Out}{\text{\rm  Out}}

\title[Structure of solutions of exponential equations]
{Structure of solutions of exponential equations in acylindrically hyperbolic groups}

\author{Agnieszka Bier}
\address{Department of Applied Mathematics, Silesian Univesity of Technology,
ul. Kaszubska 23, 44 - 101 Gliwice, Poland}
\email{agnieszka.bier@polsl.pl}

\author{Oleg Bogopolski}
\address{Institute of Mathematics, University of Szczecin, ul. Wielkopolska 15, 70-451 Szczecin, Poland}
\email{oleg.bogopolskiy@usz.edu.pl}

\begin{document}

\keywords{exponential equations, acylindrically hyperbolic groups, relatively hyperbolic groups, hyperbolically embedded subgroups, elliptic and loxodromic elements, definablity, (weak) Presburger arithmetic.}
\subjclass[2010]{Primary 20F65, 20F70; Secondary 20F67.}

\maketitle

\begin{abstract} Let $G$ be a group acting acylindrically on a hyperbolic space and let $E$ be an exponential equation over $G$.
We show that $E$ is equivalent to a finite disjunction of finite systems of pairwise
independent equations  which are either loxodromic over virtually cyclic subgroups or elliptic.
We also obtain a description of the solution set of $E$. We obtain stronger results in the case where $G$ is hyperbolic relative to a collection of peripheral subgroups $\{H_{\lambda}\}_{\lambda\in \Lambda}$. In particular, we prove in this case that the solution sets of exponential equations over $G$ are $\mathbb{Z}$-semilinear if and only if the solution sets of exponential equations over every $H_{\lambda}$, $\lambda\in \Lambda$, are $\mathbb{Z}$-semilinear.
We obtain an analogous result for finite disjunctions of finite systems of exponential equations and inequations
over relatively hyperbolic groups in terms of definable sets in the weak Presburger arithmetic.
\end{abstract}

\section{Introduction}

In 2015, Myasnikov, Nikolaev and Ushakov  initiated the study of exponential equations in groups~\cite{MNU} (see Definition~\ref{Def_Exp_Eq})
which has become a topic of intensive investigations on the edge of group theory and complexity theory~\cite{Bogo_Ivanov, KLZ, Dudkin, LZ_1, LZ_2, Frenkel, GKLZ, Figelius, MT,
Lohrey_1,Lohrey_2, FLZ}.

The results obtained in~\cite{Lohrey_2, MNU} for hyperbolic groups motivated us to investigate exponential equations in the much wider class of acylindrically hyperbolic groups. This class of groups includes
non-(virtually cyclic) groups that are hyperbolic relative to proper subgroups,
many 3-manifold groups, many groups acting on trees,
non-(virtually cyclic) groups acting properly on proper CAT(0)-spaces and containing rank-one elements, non-cyclic directly indecomposable right-angled Artin groups,
all but finitely many mapping class groups,
groups of deficiency at least~2,
$\Out(F_n)$ for $n\geqslant 2$ (see \cite{Osin_1, Osin_3} for references and historical remarks).
All necessary information about acylindrically hyperbolic groups is given in Section 2 below. We refer to the manuscript~\cite{Osin_0} for definitions concerning relatively hyperbolic groups.

\begin{definition}\label{Def_Exp_Eq}
An {\it exponential equation} over a group $G$ is an equation of the form
\begin{equation}\label{eq1}
a_1g_1^{x_1}a_2g_2^{x_2}\dots a_ng_n^{x_n}=1,
\end{equation}
where $a_1,g_1,\dots, a_n,g_n$ are elements from $G$ (called {\it coefficients}) and $x_1,\dots,x_n$ are {\it variables} (which take values in $\mathbb{Z}$).
A tuple $(k_1,\dots,k_n)$ of integers is called a {\it solution} of this equation if $a_1g_1^{k_1}a_2g_2^{k_2}\dots a_ng_n^{k_n}=1$ in $G$.
Two exponential equations are called {\it independent} if the sets of their variables are disjoint.
\end{definition}

The first two theorems say that if $G$ is acylindrically hyperbolic (respectively, relatively hyperbolic), then
equation (1) can be reduced to a finite number of (in general) simpler equations.

In the following definition we introduce two types of simpler equations -- loxodromic and elliptic equations.
In preparation for this definition, we advise the reader to review subsections~2.3 and 2.5 below.
By default, all actions of groups on metric spaces are assumed to be isometric in this paper.
If the action is fixed, we will sometimes omit the words ``with respect to this action''.

\begin{definition}\label{Lox_Ell_Eq}
Let $G$ be a group acting on a metric space $S$.
We say that the exponential equation~\eqref{eq1} is
{\it loxodromic} (\!\!\!{\it~elliptic}) with respect to this action if all $g_i$ in this equation
are loxodromic (\!{\it elliptic}).

For concise formulations, we introduce a very special type of elliptic equation:
An exponential equation is called {\it finitary} if it has the form $ag^x=1$, where $g$ has finite order and $a\in G$.

\end{definition}






Recall that if $G$ is a group acting acylindrically on a hyperbolic space $S$ and $g\in G$ is loxodromic,
then $g$ is contained in a unique maximal virtually cyclic subgroup, denoted by $E_G(g)$, see~\cite[Lemma 6.5]{DOG}.

\medskip

\noindent
{\bf Theorem A.}
{\it Let $G$ be an a group acting acylindrically on a hyperbolic space $S$.
Any exponential equation
$
a_1g_1^{x_1}a_2g_2^{x_2}\dots a_ng_n^{x_n}=1
$
over $G$ is equivalent to a finite disjunction $\Phi$ of finite conjunctions of exponential equations,
$$
\Phi:=\overset{k}{\underset{i=1}{\bigvee}}\overset{s_i}{\underset{j=1}{\bigwedge}} L_{i,j},
$$
such that equations in each conjunction are pairwise independent and each $L_{i,j}$ is either loxodromic or elliptic.
Every loxodromic $L_{i,j}$ is an exponential equation over the virtually cyclic group $E_G(g_{\ell})$ for some
loxodromic element $g_\ell\in \{g_1,\dots,g_n\}$.




}

\medskip

In the important special case where $G$ is relatively hyperbolic, one can take the next step and
reduce  elliptic exponential equations over $G$ to several systems consisting of equations over peripheral subgroups and of finitary equations over~$G$.

\medskip

\noindent
{\bf Theorem B.}
{\it Let $G$ be a group which is relatively hyperbolic with respect to a collection of its subgroups $\{H_\lambda\}_{\lambda\in \Lambda}$. Any exponential equation
$
a_1g_1^{x_1}a_2g_2^{x_2}\dots a_ng_n^{x_n}=1
$
over $G$ is equivalent to a finite disjunction $\Phi$ of finite conjunctions of exponential equations,
$$
\Phi:=\overset{k}{\underset{i=1}{\bigvee}}\overset{s_i}{\underset{j=1}{\bigwedge}} L_{i,j},
$$
such that equations in each conjunction are pairwise independent and each $L_{i,j}$ is either
over the virtually cyclic group $E_G(g_{\ell})$ associated with some loxodromic element $g_\ell\in \{g_1,\dots,g_n\}$, or
over the peripheral subgroup $H_{\lambda}$ for some $\lambda\in \Lambda$, or finitary over $G$.
}

\medskip

These theorems help to describe the structure of the solution set of an exponential equation over an acylindrically hyperbolic (respectively, over a relatively hyperbolic) group, see Corollaries~C, D, and~E below. Before we formulate these corollaries, we recall some definitions that have sources in model theory.

According to~\cite{GS}, a subset $S\subseteq \mathbb{N}^n$ is called {\it linear} if there exists a natural $d$ such that $S=\{A\overline {z}+\bar{b}\,|\, \bar{z}\in \mathbb{N}^d\}$ for some $n\times d$-matrix $A$ over $\mathbb{N}$ and some vector $\bar{b}$ from $\mathbb{N}^n$.\break A subset of $\mathbb{N}^n$ is called {\it semilinear} if it is a finite union of linear subsets. In~\cite{GS}, Ginsburg and Spanier proved that semilinear subsets of $\mathbb{N}^n$ are exactly those subsets of $\mathbb{N}^n$ that are definable in Presburger arithmetic,
which is a first-order theory for the structure $\langle\mathbb{N};+,0,1\rangle$ (see also~\cite{Haase}).




The notion of semilinearity naturally appears in studying of $\mathbb{N}$-solutions of exponential
equations (i.e., of solutions with all coordinate values in $\mathbb{N}$), see~\cite{FLZ, Lohrey_2}.
Since we are mostly interested in solutions 
with all coordinate values in $\mathbb{Z}$, we recall the following definition from~\cite{CF}.


\begin{definition}
A subset $S\subseteq \mathbb{Z}^n$  is called $\mathbb{Z}$-{\it linear} if there exists a natural $d$ such that
$S=\{A\overline {z}+\bar{b}\,|\, \bar{z}\in \mathbb{Z}^d\}$ for some $n\times d$-matrix $A$ over $\mathbb{Z}$
and some vector $\bar{b}$ from $\mathbb{Z}^n$.
A subset of $\mathbb{Z}^n$ is called $\mathbb{Z}$-{\it semilinear} if it is a finite union of $\mathbb{Z}$-linear subsets.
\end{definition}


Note that if $S\subseteq \mathbb{Z}^n$ is $\mathbb{Z}$-semilinear,
then $S\cap \mathbb{N}^n$ is semilinear in $\mathbb{N}^n$.
In~\cite[Theorem 3.1]{CF}, Chouffrut and Frigeri proved that finite unions of differences of $\mathbb{Z}$-semilinear subsets of $\mathbb{Z}^n$
are exactly those subsets of $\mathbb{Z}^n$ that are definable in the first order theory for the structure $\langle\mathbb{Z};+,0,1\rangle$ which is called  the
{\it weak Presburger arithmetic}.  The word weak comes from the fact that the predicate $<$ is not expressible by first-order formulas in this arithmetic; however it is expressible in Presburger arithmetic.





We identify the direct sum $\mathbb{Z}^s\oplus \mathbb{Z}^t$ with $\mathbb{Z}^{s+t}$ 
in the usual way, $s,t\in \mathbb{N}\cup \{0\}$.


\medskip

\noindent
{\bf Corollary C.} {\it Let $G$ be a group and $a_1g_1^{x_1}a_2g_2^{x_2}\dots a_ng_n^{x_n}=1$
an exponential equation over $G$.
\begin{enumerate}
\item[{\rm (a)}]
Suppose that $G$ acts acylindrically on a hyperbolic space $S$.
Let $\mathcal{L}$ (respectively, $\mathcal{E}$) be the set of all subscripts $j\in \{1,\dots,n\}$ such that $g_j$ is loxodromic (respectively, elliptic).
Then, for some finite $s$, there exist $\mathbb{Z}$-semilinear subsets
$N_1,\dots, N_s\subseteq \mathbb{Z}^{\mathcal{L}}$
and subsets $M_1,\dots,M_s\subseteq \mathbb{Z}^{\mathcal{E}}$ such that
the solution set of the above equation coincides, up to appropriate enumeration of components, with
$$
\overset{s}{\underset{i=1}{\cup}} (N_i\oplus M_i).
$$
More precisely, every $M_i$ is the solution set of a finite system of independent elliptic equations over $G$.

\medskip

\item[{\rm (b)}] Suppose that $G$ is relatively hyperbolic with respect to a collection of its subgroups $\{H_\lambda\}_{\lambda\in \Lambda}$. Then the conclusion in {\rm (a)} can be strengthened: Every $M_i$
    is the solution set of a finite system of independent equations, each of these equations is either over $H_{\lambda}$ for some $\lambda\in \Lambda$, or finitary over~$G$.
\end{enumerate}
}

\medskip

In~\cite{Lohrey_2}, Lohrey proved that the set of $\mathbb{N}$-solutions of an exponential equation~\eqref{eq1} over a hyperbolic group is a semilinear subset of $\mathbb{N}^n$. Also  Figelius, Lohrey and Zetzsche proved in~\cite{FLZ} that the semilinearity of the $\mathbb{N}$-solution sets of exponential equations is preserved under taking of amalgamated products over finite subgroups and under taking of HNN extensions over finite associated subgroups.
Corollary D generalizes this result to relatively hyperbolic groups.


\medskip

\noindent
{\bf Corollary D.} {\it Let $G$ be a group which is relatively hyperbolic with respect to a collection of its subgroups $\{H_\lambda\}_{\lambda\in \Lambda}$.
Then the solution set of any exponential equation over $G$ is $\mathbb{Z}$-semilinear if and only if for every $\lambda\in \Lambda$
the solution set of any exponential equation over $H_{\lambda}$ is $\mathbb{Z}$-semilinear.

The same is valid for solutions with components in $\mathbb{N}$; we shall only replace\break $\mathbb{Z}$-semi\-linear by semilinear.
}

\medskip

In model theory, it is more natural to consider finite disjunctions of finite systems of equations and inequations rather than
one equation. Therefore we introduce the following definition and formulate Corollary E for this case.

\begin{definition} An {\it exponential inequation} over a group $G$ is an expression as in (1) with $=$ replaced by $\neq$.
For brevity, we call any finite disjunction of finite systems of exponential equations and inequations over $G$ by a {\it generalized exponential (in)equation over $G$}.
\end{definition}

\noindent
{\bf Corollary E.}
{\it
Let $G$ be a group which is relatively hyperbolic with respect to a collection of its subgroups $\{H_\lambda\}_{\lambda\in \Lambda}$.
Then the solution set of any generalized exponential (in)equation over $G$ with $n$ variables
is a definable subset of $\mathbb{Z}^n$ in the weak Presburger arithmetic if and only if
the same holds for every $H_{\lambda}$, $\lambda\in \Lambda$.

The same is valid for solutions with components in $\mathbb{N}$; we shall only replace $\mathbb{Z}^n$ by $\mathbb{N}^n$ and omit the word weak.
}


\medskip

In a forthcoming paper, we estimate a ``minimal solution''
of equation~\eqref{eq1} over an acylindrically hyperbolic group and over a relatively hyperbolic group, see preprint~\cite{Bogo_Bier}.
We obtain there linear estimates (in terms of lengths of coefficients) of the minimal loxodromic part.
Polynomial and linear estimates in the case of hyperbolic groups
were obtained in~\cite{MNU}~and~\cite{Lohrey_2}, respectively.





We thank Aleksander Ivanov and David Bradley-Williams for their remarks concerning Presburger arithmetic.

\section{Definitions and preliminaries}

We introduce general notation and recall some relevant definitions and statements
from the papers~\cite{DOG,Osin_1}.

\subsection{General notation} All generating sets considered in this paper are assumed to be {\it symmetric}, i.e., closed under taking inverse elements.
Let $G$ be a group generated by a subset $X$. For $g\in G$ let $|g|_X$ be the length of a shortest word in $X$ representing $g$. The corresponding metric
on $G$ is denoted by ${d}_X$ (or by ${d}$ if $X$ is clear from the context); thus ${d}_X(a,b)=|a^{-1}b|_X$. The right Cayley graph of $G$ with respect to $X$ is denoted by $\Gamma(G,X)$.
By a path $p$ in the Cayley graph we mean a combinatorial path; the initial and the terminal vertices of $p$ are denoted by $p_{-}$ and $p_{+}$, respectively.
The length of $p$ is denoted by $\ell(p)$. The {\it label} of $p$ (which is a word in the alphabet $X$)
is denoted by ${\bold{ Lab}}(p)$.


A path $p$ in $\Gamma(G,X)$ is called ($\varkappa,\varepsilon)$-{\it quasi-geodesic}, where $\varkappa\geqslant 1$,
$\varepsilon\geqslant 0$, if ${d}(q_{-},q_{+})\geqslant \frac{1}{\varkappa}\ell(q)-\varepsilon$ for any subpath $q$ of $p$.

Recall that two elements $a,b\in G$ of infinite order are called {\it commensurable}, if there exist nonzero integers $s,t$ and an element $g\in G$ such that $a^s = g^{-1}b^tg$.

\subsection{Two equivalent definitions of acylindrically hyperbolic groups}

\begin{definition}\label{acyl_action} {\rm (see~\cite{Bowditch} and Introduction in~\cite{Osin_1})
An action of a group $G$ on a metric space $S$ is called
{\it acylindrical}
if for every $\varepsilon>0$ there exist $R,N>0$ such that for every two points $x,y\in S$ with $d(x,y)\geqslant R$,
there are at most $N$ elements $g\in G$ satisfying
$$
d(x,gx)\leqslant \varepsilon\hspace*{2mm}{\text{\rm and}}\hspace*{2mm} d(y,gy)\leqslant \varepsilon.
$$
}
\end{definition}

Given a generating set $X$ of a group $G$, we say that the Cayley graph $\Gamma(G,X)$ is
{\it acylindrical} if the left action of $G$ on $\Gamma(G,X)$ is acylindrical.
For Cayley graphs, the acylindricity condition can be rewritten as follows:
for every $\varepsilon>0$ there exist $R,N>0$ such that for any $g\in G$ of length $|g|_X\geqslant R$
we have
$$
\bigl|\{f\in G\,|\, |f|_X\leqslant \varepsilon,\hspace*{2mm} |g^{-1}fg|_X\leqslant \varepsilon \}\bigr|\leqslant N.
$$

Recall that a geodesic metric space $\frak{X}$ is called {\it $\delta$-hyperbolic} if each side of any geodesic triangle $\Delta$ in $\frak{X}$ lies in the $\delta$-neighborhood of the union of the other two sides of $\Delta$.
An action of a group $G$ on a hyperbolic space $S$ is called {\it elementary} if the limit set
of $G$ on the Gromov boundary $\partial S$ contains at most 2 points.

\begin{definition}\label{Definition_of_Osin} {\rm (see~\cite[Definition 1.3]{Osin_1})
A group $G$ is called {\it acylindrically hyperbolic} if it satisfies one of the following equivalent
conditions:

\begin{enumerate}
\item[(${\rm AH}_1$)] There exists a generating set $X$ of $G$ such that the corresponding Cayley graph $\Gamma(G,X)$
is hyperbolic, $|\partial \Gamma (G,X)|>2$, and the natural action of $G$ on $\Gamma(G,X)$ is acylindrical.

\medskip

\item[(${\rm AH}_2$)] $G$ admits a non-elementary acylindrical action on a hyperbolic space.
\end{enumerate}
}
\end{definition}

In the case (AH$_1$), we also write that $G$ is {\it acylindrically hyperbolic with respect to $X$}.

\medskip

\subsection{Elliptic and loxodromic elements}

The following definition is standard.

\begin{definition}\label{elliptic_and_loxodromic}
{\rm
Let $G$ be a group acting on a metric space $S$, in symbols $G\curvearrowright S$. An element $g\in G$ is called {\it elliptic} with respect to this action if some (equivalently, any) orbit of $g$ in $S$ is bounded.
An element $g\in G$ is called {\it loxodromic} with respect to this action if the map
$\mathbb{Z}\rightarrow S$ defined by
$n\mapsto g^nx$ is a quasi-isometric embedding for some (equivalently, any) $x\in S$. That is,
for $x\in S$, there exist $\varkappa\geqslant 1$ and $\varepsilon\geqslant 0$ such that for any $n,m\in \mathbb{Z}$ we have
$$
d(g^nx,g^mx)\geqslant \frac{1}{\varkappa} |n-m|-\varepsilon.
$$

The sets of all elliptic and all loxodromic elements of $G$ with respect to the action $G\curvearrowright S$ is denoted by ${\rm Ell}(G\curvearrowright S)$ and ${\rm Lox}(G\curvearrowright S)$, respectively.

Let $X$ be a generating set of $G$.
We say that $g\in G$ is {\it elliptic} (respectively, {\it loxodromic}) {\it with respect to $X$} if $g$ is elliptic (respectively, loxodromic) for the canonical left action of $G$ on the Cayley graph $\Gamma(G,X)$.
The set of all elliptic and all loxodromic elements of $G$ with respect to $X$ are denoted by ${\rm Ell}(G,X)$ and ${\rm Lox}(G,X)$, respectively.

If the action $G\curvearrowright S$ or the generating set $X$ of $G$ are clear from a context, we omit the words ``with respect to ...''.
}
\end{definition}

Note that for groups acting on geodesic hyperbolic spaces,
there is only one additional isometry type of an element- parabolic
(see e.g.~\cite[Chapitre 9, Th$\acute{\text{e}}$or$\grave{\text{e}}$me~2.1]{CDP}).



Now suppose that $G$ is a group acting acylindrically on a hyperbolic space $S$. Bowditch~\cite[Lemma 2.2]{Bowditch} proved that in this situation every element $g$ of $G$ is either elliptic or loxodromic (see a more general statement in~\cite[Theorem 1.1]{Osin_1}).
If $g$ is loxodromic, then it is contained in a
unique maximal virtually cyclic subgroup~\cite[Lemma 6.5]{DOG}.
This subgroup, denoted by $E_G(g)$, is called the {\it elementary subgroup associated with $g$}; it can be described as follows (see equivalent definitions in~\cite[Corollary~6.6]{DOG}):
$$
\hspace*{12.5mm}\begin{array}{ll}
E_G(g)\! \!\! & =\{f\in G\,|\, \exists\,  n\in \mathbb{N}:  f^{-1}g^nf=g^{\pm n}\}.
\vspace*{3mm}\\
\! \!\! & =\{f\in G\,|\, \exists\,  k,m\in \mathbb{Z}\setminus \{0\}:  f^{-1}g^kf=g^{m}\}.
\end{array}
$$




\subsection{Weakly hyperbolic groups}
Let $G$ be a group, $\{H_{\lambda}\}_{\lambda\in \Lambda}$ a collection of subgroups of $G$.
A subset $X$ of $G$ is called a {\it relative generating set of $G$ with respect to}
$\{H_{\lambda}\}_{\lambda\in \Lambda}$ if $G$ is generated by $X$ together with the union of all $H_{\lambda}$.
In this paper we assume that all relative generating sets are {\it symmetric}, i.e., for every $x\in X$, we have $x^{-1}\in X$. For technical reasons we also assume that $1\in X$.
We define
$$
\mathcal{H}=\bigsqcup_{\lambda\in\Lambda}H_{\lambda}.
$$


\noindent
{\bf Assumption.} In Definitions~\ref{one} --~\ref{n-gon} we assume that $X$ is a relative generating set of $G$ with respect to $\{H_{\lambda}\}_{\lambda\in \Lambda}$.

\begin{definition}\label{one} {\rm (see~\cite[Definition 4.1]{DOG})
The group $G$ is called {\it weakly hyperbolic} relative to $X$ and $\{H_{\lambda}\}_{\lambda\in \Lambda}$ if the Cayley graph $\Gamma(G, X\sqcup \mathcal{H})$ is hyperbolic.
}
\end{definition}

\noindent
We consider the Cayley graph
$\Gamma(H_{\lambda},H_{\lambda})$ as a complete subgraph of $\Gamma(G,X\sqcup \mathcal{H})$.

\begin{definition}\label{definition_relative_metrics}
{\rm (see~\cite[Definition 4.2]{DOG})
For every $\lambda\in \Lambda$, we introduce a {\it relative metric}
$\widehat{d}_{\lambda}:H_{\lambda}\times H_{\lambda}\rightarrow [0,+\infty]$ as follows:

Let $a,b\in H_{\lambda}$. A path
in $\Gamma(G,X\sqcup \mathcal{H})$ from $a$ to $b$ is called {\it $H_{\lambda}$-admissible} if it has no edges in the subgraph $\Gamma(H_{\lambda},H_{\lambda})$.

The distance $\widehat{d}_{\lambda}(a,b)$ is defined to be the length of a shortest
{\it $H_{\lambda}$-admissible} path connecting $a$ to $b$ if such exists.
If no such path exists, we set $\widehat{d}_{\lambda}(a,b)=\!
\infty$.


}
\end{definition}




\begin{definition}\label{components}
{\rm (see~\cite[Definition 4.5]{DOG})
Let $q$ be a path in the Cayley graph $\Gamma(G,X\sqcup \mathcal{H})$. A non-trivial subpath $p$ of $q$
is called an {\it $H_{\lambda}$-subpath} if the label of $p$ is a word in the alphabet $H_{\lambda}$.
An $H_{\lambda}$-subpath $p$ of $q$ is called an {\it $H_{\lambda}$-component} if $p$ is not contained in a longer subpath of $q$ with this property. Further by a {\it component} of $q$ we mean an $H_{\lambda}$-component of $q$ for some $\lambda\in \Lambda$. Two $H_{\lambda}$-components $p_1,p_2$ of a path $q$ in $\Gamma(G,X\sqcup \mathcal{H})$ are called {\it connected} if there exists a path $\gamma$ in $\Gamma(G,X\sqcup \mathcal{H})$ that connects some vertex of $p_1$ to some vertex of $p_2$, and ${\bf Lab}(\gamma)$ is a word consisting only of letters from
$H_{\lambda}$.

Note that we can always assume that $\gamma$ has length at most 1 as every element of $H_{\lambda}$
is included in the set of generators. An $H_{\lambda}$-component $p$ of a path $q$ in $\Gamma(G,X\sqcup \mathcal{H})$ is {\it isolated} if it is not connected to any other component of $q$.

}
\end{definition}

Given a path $p$ in $\Gamma(G, X\sqcup \mathcal{H})$, the canonical image of ${\bold{Lab}}(p)$ in $G$ is denoted by ${\bold {Lab}}_G(p)$.
By $X^{\ast}$ we denote the free monoid generated by $X$.



\begin{definition}\label{n-gon} {\rm (see~\cite[Definition 4.13]{DOG})
Let $\varkappa \geqslant 1$, $\varepsilon\geqslant 0$, and $m\geqslant 2$. Let $\mathcal{P}=p_1\dots p_m$ be an $m$-gon in $\Gamma(G, X\sqcup \mathcal{H})$ and let $I$ be a subset of the set of its sides $\{p_1,\dots,p_m\}$
such that:

1) Each side $p_i\in I$ is an isolated $H_{\lambda_i}$-component of $\mathcal{P}$ for some $\lambda_i\in \Lambda$.

2) Each side $p_i\notin I$ is a $(\varkappa,\varepsilon)$-quasi-geodesic.

\medskip
\noindent
We denote $s(\mathcal{P},I)=\underset{p_i\in I}{\sum} \widehat{d}_{\lambda_i}(1,{\bold{Lab}}_G(p_i))$.

}
\end{definition}

\begin{proposition}\label{Proposition_Osin}
{\rm (see~\cite[Proposition 4.14]{DOG})}
Suppose that $G$ is weakly hyperbolic relative to $X$ and $\{H_{\lambda}\}_{\lambda\in \Lambda}$.
Then for any $\varkappa\geqslant 1$, $\varepsilon\geqslant 0$, there exists a constant $D(\varkappa,\varepsilon)>0$ such that
for any $m$-gon $\mathcal{P}$ in $\Gamma(G,X\sqcup \mathcal{H})$ and any subset $I$ of the set of its sides satisfying conditions
of Definition~\ref{n-gon}, we have $s(\mathcal{P},I)\leqslant D(\varkappa,\varepsilon)m$.
\end{proposition}

The following definition is a stronger version of the weak hyperbolicity. It will be used in the proof
of Theorem~A to provide there some finiteness conditions.

\begin{definition}\label{def_hyperb_embedd} {\rm (see~\cite[Definition 4.25]{DOG})
Let $G$ be a group, $X$ a (not necessarily finite) subset of $G$. A collection of subgroups $\{H_{\lambda}\}_{\lambda\in \Lambda}$
of $G$ is called {\it hyperbolically embedded in $G$ with respect to $X$}
(we write $\{H_{\lambda}\}_{\lambda\in \Lambda}\hookrightarrow_h (G,X)$) if the following hold.

\begin{enumerate}
\item[(a)] The group $G$ is generated by $X$ together with the union of all $H_{\lambda}$ and the Cayley graph
$\Gamma(G,X\sqcup \mathcal{H})$ is hyperbolic.

\item[(b)] For every $\lambda\in \Lambda$, the metric space $(H_{\lambda},\widehat{d}_{\lambda})$ is
proper. That is, any ball of finite radius in $H_{\lambda}$ contains finitely many elements.
\end{enumerate}
}
\end{definition}

\medskip

By~\cite[Theorem 1.2]{Osin_1}, a group $G$ is acylindrically hyperbolic if and only if $G$ contains a proper infinite hyperbolically embedded subgroup.

\newpage

\subsection{Loxodromic and elliptic elements in relatively hyperbolic groups.}
Suppose that $G$ is relatively hyperbolic with respect to a collection of its subgroups $\{H_{\lambda}\}_{\lambda\in \Lambda}$ (see Definition~2.35 in~\cite{Osin_0}). In particular, there exists a finite set $X$ such that $G$ is generated by the set
$$
Y=X\cup (\underset{\lambda\in \Lambda}{\cup} H_{\lambda}).
$$

\noindent
It is known that any element of $G$ is either elliptic or loxodromic with respect to~$Y$. Moreover, an element of $G$ is elliptic with respect to $Y$ if and only if it is conjugate to an element of one of the subgroups $H_{\lambda}$, $\lambda\in \Lambda$, or has a finite order (see Definition~4.1 and Theorem~4.25 in~\cite{Osin_0}). Clearly, this classification does not depend on a choice of a finite relative generating set $X$.

\section{Proof of Theorems A and B}
\begin{definition} (Catalan permutations)
A permutation $\sigma$ of the set $\{1,\dots,n\}$ is called a {\it Catalan permutation} if the orbits of $\sigma$ determine a non-crossing partition of this set. We recall that two blocks $\mathcal B _1, \mathcal B_2$ of a partition of $\{1,\dots,n\}$ are called {\it crossing} if there exist $k,l\in \mathcal B_1$ and $s,t\in \mathcal B_2$ such that either $k<s<l<t$ or $s<k<t<l$. The set of all such permutations is denoted by $\Cat(n)$.
\end{definition}

It is known that the number of Catalan permutations of $\{1,\dots,n\}$ is equal to the $n$-th Catalan number $$
C_n=\frac{1}{n+1}\binom{2n}{n}.
$$

\begin{lemma}\label{Catalan_connection}
Let $G$ be a group, $\{H_{\lambda}\}_{\lambda\in \Lambda}$ a collection of subgroups of $G$, and $X$ a symmetric relative generating set of $G$ with respect to $\{H_{\lambda}\}_{\lambda\in \Lambda}$ with $1\in X$.\break
Let $\mathcal{P}=p_1e_1p_2e_2\dots p_ne_n$ be a cycle in the Cayley graph $\Gamma=\Gamma(G,X\sqcup \mathcal{H})$,
where $\mathcal{H}=\underset{\lambda\in \Lambda}{\sqcup} H_{\lambda}$,
such that each $p_i$ is a nontrivial path in $\Gamma$ with $\Lab(p_i)\in X^{\ast}$ and each~$e_i$ is an edge in $\Gamma$ with $\Lab(e_i)\in H_{\lambda(i)}$ for some $\lambda(i)\in \Lambda$, $i=1,\dots,n$. Then there exist a Catalan permutation $\sigma$ of $\{1,\dots,n\}$
and edges $f_1,\dots ,f_n$ in $\Gamma$ such that for any $i$ the following holds:


\begin{enumerate}
\item[{\rm (1)}] $f_i$ connects $(e_i)_{+}$ with $(e_{\sigma(i)})_{-}$, and $\Lab(f_i)\in H_{\lambda(i)}$.
\item[{\rm (2)}]
If $\ell$ is the length of the $\sigma$-orbit of $i$, then the label of every  edge in the cycle $e_if_ie_{\sigma(i)}f_{\sigma(i)}\dots e_{\sigma^{\ell-1}(i)}f_{\sigma^{\ell-1}(i)}$ lies in
$H_{\lambda(i)}$.
\item[{\rm (3)}] $f_i$ is an isolated $H_{\lambda(i)}$-component of some cycle in $\Gamma$ of length at most the length of $\mathcal{P}$.

\end{enumerate}
\end{lemma}

\vspace*{-15mm}
\hspace*{15mm}
\includegraphics[scale=0.5]{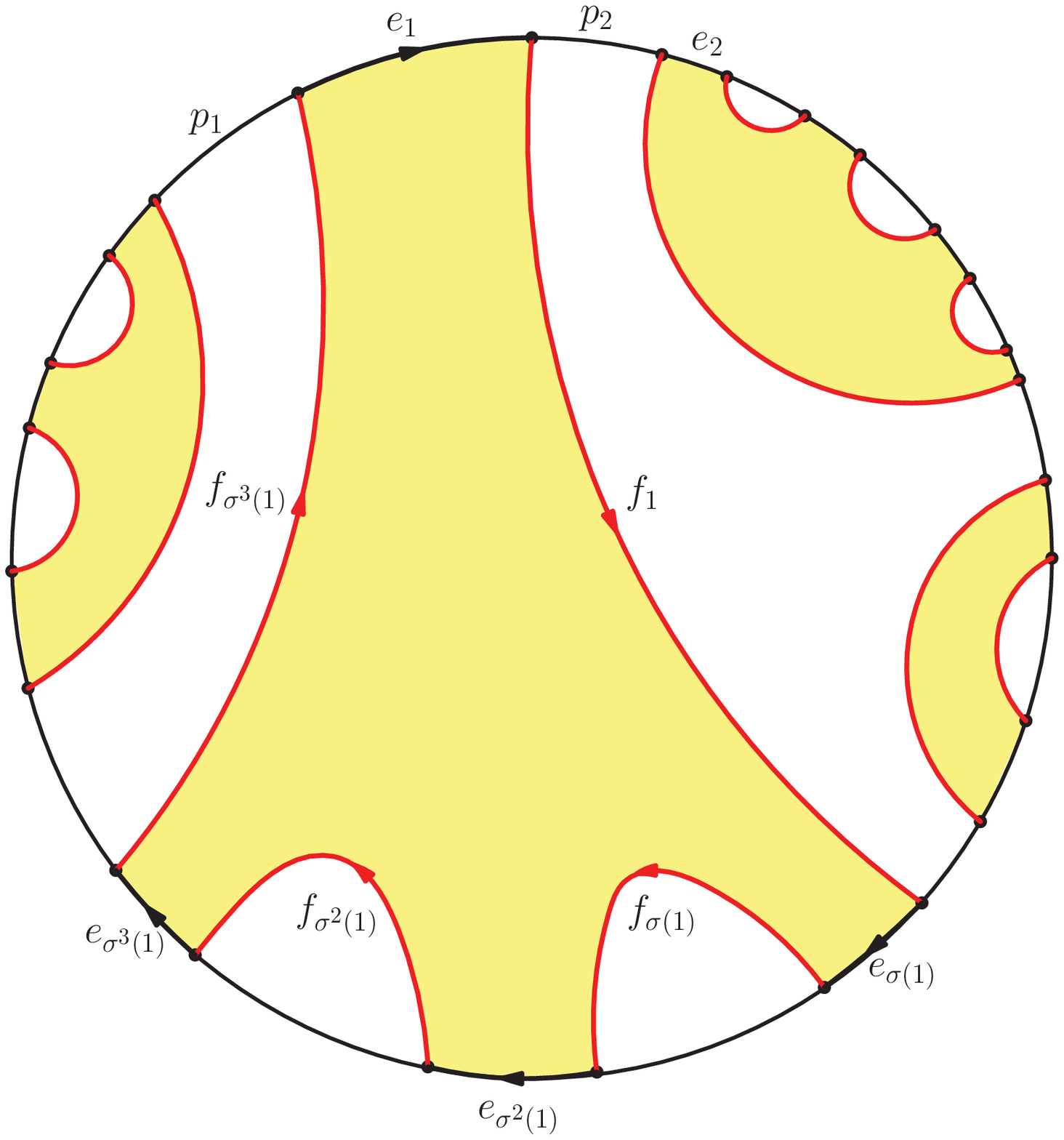}

\vspace*{-35mm}

\begin{center}
Fig. 1. Illustration to Lemma~\ref{Catalan_connection}.
\end{center}

\begin{remark}
Condition (2) implies that if $\ell>1$, then $e_i$ is connected with every edge from the set $\{e_{\sigma(i)}, \dots, e_{\sigma^{\ell-1}(i)}\}$. Therefore if $e_i$ is isolated in $\mathcal{P}$, then $\ell=1$, i.e., $i$ is a fixpoint of $\sigma$.
\end{remark}

{\it Proof of lemma.} Induction on $n$.
Suppose that $n=1$. Then the cycle $\mathcal P = p_1e_1$ contains only one $H_{\lambda}$-component, namely $e_1$, which is necessarily isolated in $\mathcal{P}$.\break In this case we set $\sigma = id$ and $f_1 = e_1^{-1}$.

Now, assume that $n>1$ and the statement is valid for all cycles with less than~$n$ components. Consider the cycle $\mathcal{P}=p_1e_1p_2e_2\dots p_ne_n$ defined in the lemma.

If all components in $\mathcal P$ are isolated, then the statement trivially holds with $\sigma = id$ and $f_i = e_i^{-1}$ for $i=1,2,\dots ,n$.
Now suppose that $\mathcal{P}$ contains a non-isolated $H_\lambda$-component $e_{i_1}$ for some $\lambda\in \Lambda$.
Let $e_{i_2},\dots,e_{i_k}$ be all $H_{\lambda}$-components of $\mathcal{P}$ which are connected with $e_{i_1}$.
We assume that $i_1<i_2<\dots <i_k<i_1$ with respect to the natural
cyclic order on $\{1,2,\dots,n\}$.
We set $\mathcal{I}=\{i_1,i_2,\dots,i_k\}$ and let $\tau=(i_1,i_2,\dots,i_k)$ be the corresponding cyclic permutation.

Then all these $H_\lambda$-components $e_{i_1}$, $e_{i_2}$,\dots , $e_{i_k}$ are pairwise connected, hence one can choose edges $f_{i_1}$, $f_{i_2}$,\dots, $f_{i_k}$ in $\Gamma$ with labels in $H_\lambda$
such that $f_{i_j}$ connects $(e_{i_j})_{+}$ with $(e_{i_{j+1}})_{-}$ for $j=1,2,\dots,k$ (addition modulo $k$).
This way we obtain a cycle in $\Gamma$ whose all edges are labelled by elements from $H_{\lambda}$:
$$\mathcal{C}_{\mathcal \tau} = e_{i_1}f_{i_1}e_{\tau (i_1)}f_{\tau (i_1)}\dots e_{\tau^{k-1} (i_1)}
f_{\tau^{k-1} (i_1)}.$$

The subset $\mathcal I$ defines a partition $\Pi$ on the set $\{1,2,\dots ,n\}\setminus \mathcal I$ into $l\leqslant k$ disjoint subsets (as some of the following sets may vanish):
$$\{i_1+1, i_1+2,\dots ,i_2-1\}\sqcup \dots \sqcup \{i_k+1, i_k+2,\dots ,n,1,\dots ,i_1-1\}$$
and $l$ respective cycles in $\Gamma$ (only for nonempty partition blocks as above):
$$\begin{array}{rcl}
\mathcal{P}_{i_1} &=& p_{i_1+1}e_{i_1+1}p_{i_1+2}e_{i_1+2} \dots p_{i_2-1}f_{i_1}^{-1}\\
\mathcal{P}_{i_2} &=& p_{i_2+1}e_{i_2+1}p_{i_2+2}e_{i_2+2} \dots  p_{i_3-1}f_{i_2}^{-1}\\
\vdots & &\\
\mathcal{P}_{i_k} &=& p_{i_k+1}e_{i_k+1}p_{i_k+2}e_{i_k+2} \dots p_{i_1-1}f_{i_k}^{-1} \end{array}$$
Note that since $e_{i_2},\dots , e_{i_k}$ are all $H_\lambda$-components connected to $e_{i_1}$, the connecting edge $f_{i_j}$ is necessarily an isolated component of $\mathcal P_{i_j}$. Moreover, $\mathcal P_{i_j}$ is strictly shorter than $\mathcal P$, as it is formed from the edge $f_{i_j}$ and a proper subpath of $\mathcal P$ which misses at least two connected $H_\lambda$-components in $\mathcal P$.
By inductive assumption, for every cycle $\mathcal P_{i_j}$, there exist a Catalan permutation $\sigma_{i_j}$ of the block of indices $\{i_{j}+1, i_j+2, \dots ,i_{j+1}-1\}$ and there exist edges  $f_{i_j+1}$, $f_{i_j+2}$, \dots, $f_{i_{j+1}-1}$ which, together with $f_{i_j}$ satisfy the statements of the lemma.

Now, let us define a permutation $\sigma$ of $\{1,2,\dots ,n\}$ as follows:
$$\sigma (s) = \begin{cases}
\tau (s),\qquad s\in \mathcal I\\
\sigma_{i_j}(s),\qquad s\in  \{i_{j}+1, i_j+2,\dots ,i_{j+1}-1\}, \ j\in 1,2,\dots ,k\\
\end{cases} $$
Then $\sigma$ is a Catalan permutation since it induces a Catalan permutation on each block and these blocks are pairwise non-crossing. This permutation and the constructed edges $f_1,\dots ,f_n$ satisfy conditions (1)-(3). \hfill $\Box$


\begin{definition}\label{def_ass_sys} (Associated systems)
Let $G$ be a group and
$$
a_1g_1^{x_1}a_2g_2^{x_2}\dots a_ng_n^{x_n}=1
$$
an exponential equation of form \eqref{eq1} over $G$.
Suppose that $\sigma$ is a Catalan permutation of the set $I=\{1,\dots,n\}$ and
$b=(b_1,\dots,b_n)$ is an $n$-tuple of elements of $G$.
We define a $(\sigma,b)$-system of equations associated with \eqref{eq1} as follows.

Take a closed disc $D$ and divide its oriented boundary $\partial D$ into $2n$ oriented segments (edges): $\partial D=e_1'e_1\dots e_n'e_n$. Label each edge $e_i'$ by $a_i$ and each edge $e_i$ by $g_i^{x_i}$.
Thus, the label of $\partial D$ is exactly the expression on the left side of  \eqref{eq1}.
For any $i\in I$, we draw an oriented chord $f_i$ from $(e_i)_{+}$ to $(e_{\sigma(i)})_{-}$
and label it by the element $b_i$. These chords divide the disc $D$ into regions, say $D_1,\dots ,D_s$ (see Fig.~2). The system of equations
$$
\Phi_{(\sigma,b)}:\overset{s}{\underset{i=1}{\bigwedge}} (\Lab(\partial D_i)=1)
$$
is called the $(\sigma,b)$-{\it system associated} with the equation \eqref{eq1}.
\end{definition}

\vspace*{-20mm}
\hspace*{10mm}
\includegraphics[scale=0.5]{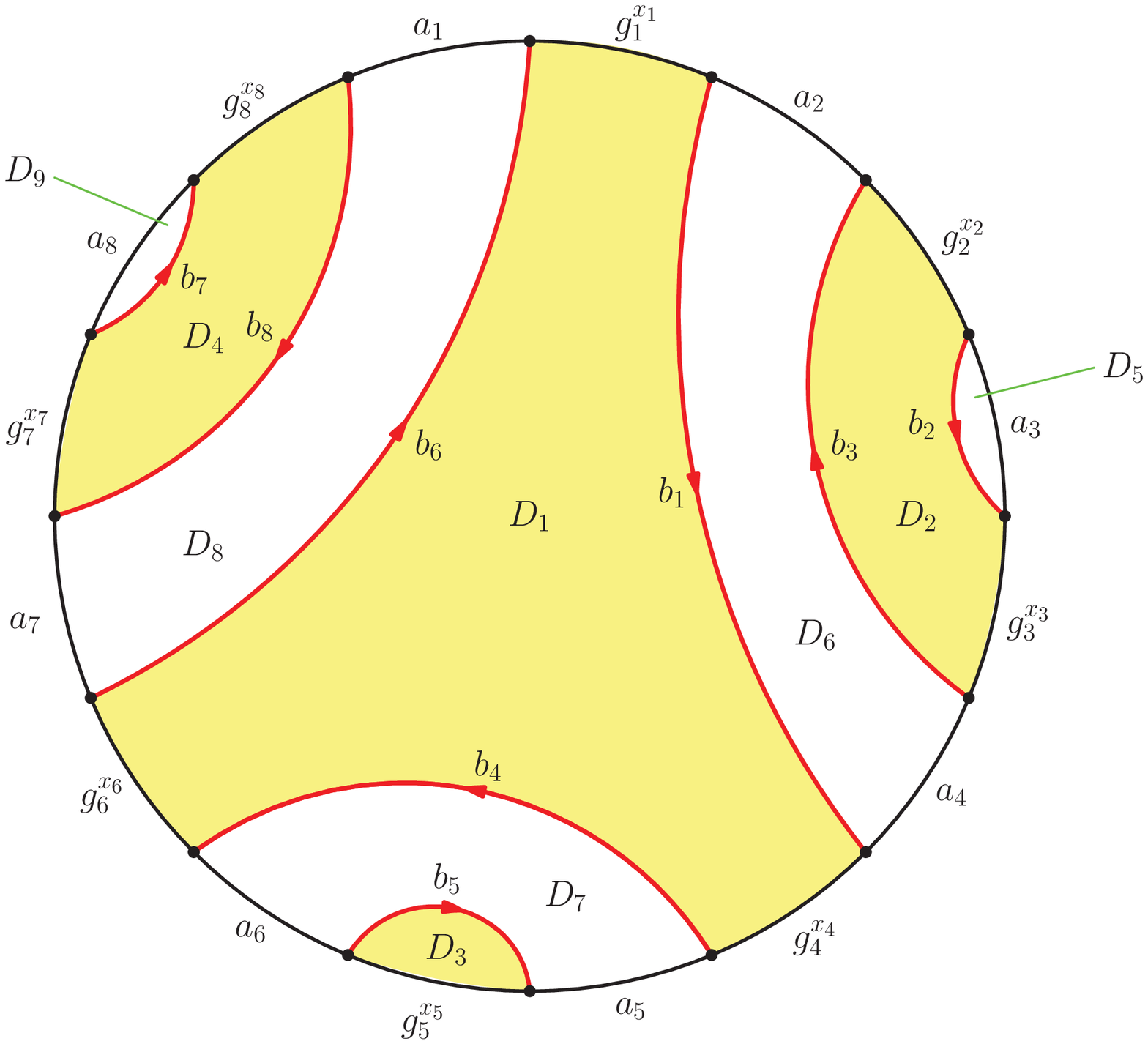}

\vspace*{-40mm}

\begin{center}
Fig. 2. Illustration to Definition~\ref{def_ass_sys}.
\end{center}

Clearly, equation (1) over an (infinite) group $G$ is equivalent to the (infinite) disjunction
$$
\underset{(\sigma,b)\in \Cat(n)\times G^n}{\bigvee}\Phi_{(\sigma,b)}.
$$
The following proposition says that if the group $G$ and equation (1) satisfy some natural conditions,
then this disjunction can be replaced by a finite one, where each system $\Phi_{(\sigma,b)}$ consists of exponential equations over some proper subgroups of $G$.






\begin{proposition}\label{bounded_coefficients}
Let $G$ be a group, $\{H_{\lambda}\}_{\lambda\in \Lambda}$ a collection of subgroups of~$G$, and $X$ a symmetric relative
generating set of $G$ with respect to $\{H_{\lambda}\}_{\lambda\in \Lambda}$ with $1\in X$.
Suppose that $\{H_{\lambda}\}_{\lambda\in \Lambda}$ is hyperbolically embedded in $G$ with respect to $X$.

Then for any two natural numbers $n$ and $C$, there exists a collection of finite sets $\{H'_{\lambda}\}_{\lambda\in \Lambda}$, where $H'_{\lambda}\subseteq H_{\lambda}$,
such that any exponential equation of form \eqref{eq1}
$$
a_1g_1^{x_1}a_2g_2^{x_2}\dots a_ng_n^{x_n}=1
$$
\noindent
with $a_i\in \langle X\rangle$, $|a_i|_X<C$, and $g_i\in H_{\lambda(i)}$ for some $\lambda(i) \in \Lambda$, $i=1,\dots,n$, is equivalent to
the finite disjunction of finite systems of independent equations
$$
\underset{(\sigma,b)\in \frak{A}}{\bigvee}\Phi_{(\sigma,b)},
$$
where
$\frak{A}=\Cat(n)\times (H'_{\lambda(1)}\times \dots \times H'_{\lambda(n)})$. Every equation in the system $\Phi_{(\sigma,b)}$ has one of the following two forms:

\medskip

$(\ast)\hspace*{37mm} g_{i_1}^{x_{i_1}}b_{i_1}g_{i_2}^{x_{i_2}}b_{i_2}\dots g_{i_s}^{x_{i_s}}b_{i_s}=1,$

\medskip
\noindent
where $1\leqslant i_1<i_2<\dots< i_s\leqslant n$ and $\lambda(i_1)=\dots =\lambda(i_s)$
and $b_{i_1},b_{i_2},\dots ,b_{i_s}\in H'_{\lambda(i_1)}$; in particular, this equation is over $H_{\lambda (i_1)}$.

\medskip

$(\ast\ast)\hspace*{36mm} a_{j_1}b^{-1}_{\ell_1}a_{j_2}b^{-1}_{\ell_2}\dots a_{j_t}b^{-1}_{\ell_t}=1,$

\medskip
\noindent
where $1\leqslant j_1<j_2<\dots< j_t\leqslant n$ and $b_{\ell_1},b_{\ell_2},\dots ,b_{\ell_t}\in  \overset{n}{\underset{i=1}{\bigcup}} H'_{\lambda(i)}$; in particular, this equation does not contain variables $x_1,\dots,x_n$.

\end{proposition}

\medskip

{\it Proof.} Let $D=D(1,1)$ be the constant from Proposition~\ref{Proposition_Osin}.
Note that this constant depends only on $G$, $X$ and $\{H_{\lambda}\}_{\lambda\in \Lambda}$.
We show that we can take
$$
H'_{\lambda}=\{h\in H_{\lambda}\,|\, \widehat{d}_{\lambda}(1, h)\leqslant Dn(C+1)\}.
$$
Note that $H'_{\lambda}$ is finite, since the metric space $(H_{\lambda}, \widehat{d}_{\lambda})$ is locally finite.
We denote by $\Gamma$ the Cayley graph $\Gamma(G,X\sqcup \mathcal{H})$, where $\mathcal{H}=\underset{\lambda\in \Lambda}{\sqcup} H_{\lambda}$.

Suppose that $k=(k_1,\dots,k_n)$ is a solution of the above equation.
We shall prove that $k$ is a solution
of an appropriate system of equations $\Phi_{(\sigma,b)}$ with $(\sigma,b)\in \frak{A}$. We represent the identity $a_1g_1^{k_1}\dots a_ng_n^{k_n}=1$ by a cycle in
the Cayley graph $\Gamma$ as follows.
For every element $a_i$ we choose a word $A_i$ of minimal positive length in the alphabet $X$ representing $a_i$ (if $a_i=1$, then $|A_i|_X=1$ since $1\in X$); in particular, $|A_i|_X\leqslant C$.
In $\Gamma$, we consider a cycle $\mathcal{P}=p_1e_1\dots p_ne_n$ such that
$p_i$ is a path with the label $A_i$ and $e_i$ is an edge with the label $g_i^{k_i}\in H_{\lambda(i)}$.
Then $\mathcal{P}$ is an $m$-gon, where
$$
m=n+\overset{n}{\underset{i=1}{\sum}}|A_i|_X\leqslant n(C+1),
$$
and the label of $\mathcal{P}$ represents the trivial element $a_1g_1^{k_1}\dots a_ng_n^{k_n}$.
We may assume that the paths $p_i$ are nontrivial, since at each vertex of $\Gamma$ there is a loop with label $1\in X$.
Then every $e_i$ is an $H_{\lambda(i)}$-component of $\mathcal{P}$.

By Lemma~\ref{Catalan_connection} applied to the cycle $\mathcal{P}$, there exist a Catalan permutation $\sigma$ of $\{1,\dots,n\}$ and edges $f_1,\dots ,f_n$ in $\Gamma$ satisfying properties (1) -- (3) of this lemma.
In particular, $f_i$ is an isolated $H_{\lambda(i)}$-component in some cycle of $\Gamma$ of length at most~$m$.
We denote $b_i=\Lab(f_i)$ and $b=(b_1,\dots,b_n)$.
By Proposition~\ref{Proposition_Osin},
$$
\widehat{d}_{\lambda(i)}(1, b_i)\leqslant Dm\leqslant Dn(C+1).
$$
Therefore $b_i\in H'_{\lambda(i)}$, and hence $(\sigma,b)\in \frak{A}$. We conclude that $(k_1,\dots,k_n)$ is a solution of the system $\Phi_{(\sigma,b)}$ and every equation in $\Phi_{(\sigma,b)}$ has the form $(\ast)$ or $(\ast\ast)$ (see Fig. 2).
\hfill $\Box$



\begin{lemma}\label{Lemma_8.8} Let $G$ be a group and $g_1,\dots,g_k\in G$ a collection of
elements of infinite order. Let $E_i$ be a virtually cyclic subgroup containing $g_i$, $i=1,\dots,k$. Suppose that
$\{E_1,\dots,E_k\}\hookrightarrow_h (G,Y)$ for some relative generating set $Y$.

Then there exists a generating set $X$ of $G$ such that $\{E_1,\dots,E_k\}\hookrightarrow_h (G,X)$,
the Cayley graph $\Gamma=\Gamma(G,X)$ is hyperbolic, the action $G\curvearrowright \Gamma$ is acylindrical,
and all $g_i$  act loxodromically on $\Gamma$.
\end{lemma}

{\it Proof.} We set $\mathcal{E}=E_1\sqcup \dots \sqcup E_k$.
By~\cite[Theorem 5.4]{Osin_1}, there exists a subset $Z\subseteq G$ such that $Y\subseteq Z$ and the following hold:

\begin{enumerate}
\item[(a)] $\{E_1,\dots ,E_k\}\hookrightarrow_h(G,Z)$. In particular, the Cayley graph $\Gamma(G, Z\sqcup \mathcal{E})$ is hyperbolic.

\item[(b)] The action of $G$ on $\Gamma(G, Z\sqcup \mathcal{E})$ is acylindrical.
\end{enumerate}

For every $i=1,\dots, k$, we choose a finite generating set $Z_i$ of $E_i$.
Then the set $X=Z\cup Z_1\cup \dots \cup Z_k$ is a generating set of $G$. Since the symmetric difference $X\Delta Z$ is finite, we still have $\{E_1,\dots,E_k\}\hookrightarrow_h (G,X)$ (see Corollary~4.27 from~\cite{DOG}).

Finally, the proof of Theorem 1.4 in Section 6 of~\cite{Osin_1} shows
that the Cayley graph $\Gamma=\Gamma(G,X)$ is hyperbolic, the action $G\curvearrowright \Gamma$ is acylindrical,
and all $g_i$  act loxodromically on $\Gamma$.
\hfill $\Box$

\medskip



{\it Proof of Theorem A.} Let $\mathcal{L}$ be the set of all indices $j\in \{1,\dots,n\}$ such that $g_j$ is loxodromic with respect to the action $G\curvearrowright S$.
We may assume $\mathcal{L}\neq \emptyset$.

{\it Step 1.} (Finding an appropriate generating set of $G$)\\
Let $\mathcal{M}$ be a maximal subset of $\{1,\dots ,n\}$ containing $\mathcal{L}$ such that
there exists an acylindrical action of $G$ on a hyperbolic  metric space for which all $g_i$, $i\in \mathcal{M}$,  act loxodromically.





Since $\mathcal{L}$ satisfies these properties and is nonempty, the set $\mathcal{M}$ is nonempty too.
We choose a maximal subset $\mathcal{M}^{nc}\subseteq \mathcal{M}$ such that the elements $g_i$ with $i\in \mathcal{M}^{nc}$ are pairwise non-commensurable.
By~\cite[Theorem 6.8]{DOG}, we have $\{E_G(g_i)\,|\, i\in \mathcal{M}^{nc}\}\hookrightarrow_h (G,Y)$ for some relative generating set $Y$.
Then, by Lemma~\ref{Lemma_8.8}, there exists a generating set $X$ of $G$ such that the following two properties are satisfied:

\begin{enumerate}
\item[{\bf (a)}] $\{E_G(g_i)\,|\, i\in \mathcal{M}^{nc}\}\hookrightarrow_h (G,X)$,

\item[{\bf (b)}] $\Gamma=\Gamma(G,X)$ is hyperbolic, the action $G\curvearrowright \Gamma$ is acylindrical, and all $g_i$, $i\in \mathcal{M}$,  act loxodromically on $\Gamma$.
\end{enumerate}

\medskip

\noindent
By maximality of $\mathcal{M}$, we have:

\begin{enumerate}
\item[{\bf (c)}] all $g_j$ with $j\in \{1,\dots,n\}\setminus \mathcal{M}$ act elliptically on $\Gamma$.
\end{enumerate}

\medskip





\medskip

{\it Step 2.} (Appropriate conjugation of loxodromic factors)\\
For every $j\in \mathcal{M}$, there exists $i\in \mathcal{M}^{nc}$ such that $g_j$ and $g_i$ are commensurable.
Then, by definition, there exists $h_j\in G$ such that $g_i^s=h_j^{-1}g_j^th_j$ for some nonzero $s,t$.
We may assume that $h_j=1$ if $j=i$. Denote $g_j'=h_j^{-1}g_jh_j$.
Then, for $j$ and $i$ as above, the elementary subgroups associated with $g_j'$ and $g_i'$ become equal:
$$
E_G(g_j')=E_G((g_j')^t)=E_G(g_i^s)=E_G(g_i)=E_G(g_i').
$$
In particular,

\begin{enumerate}
\item[{\bf (d)}] for every $j\in \mathcal{M}$, there exists $i\in \mathcal{M}^{nc}$ such that $g_j'\in E_G(g_i)$.
\end{enumerate}


\medskip

{\it Step 3.} (Grouping of elliptic factors and applying Proposition~\ref{bounded_coefficients})\\
We call variables $x_j$ with $j\in \mathcal{M}$ {\it loxodromic}, and with $j\in \{1,\dots,n\}\setminus \mathcal{M}$ {\it elliptic}.
Let $\mathcal{M}=\{i_1,\dots ,i_s\}$, i.e., $g_{i_1},\dots,g_{i_s}$ are all elements of the set $\{g_1,\dots,g_n\}$ which act loxodromically on $\Gamma$. We assume $i_1<i_2<\dots <i_s$.
First we write equation~\eqref{eq1}
in the form
\begin{equation}\label{eq9}
u_1g_{i_1}^{x_{i_1}}\dots u_sg_{i_s}^{x_{i_s}}=1,
\end{equation}
where each $u_j$ is the subword of the left side of equation \eqref{eq1} between $g_{i_{j-1}}^{x_{i_{j-1}}}$ and $g_{i_j}^{x_{i_j}}$ (considered cyclically). Note that $u_j$ is a product of some coefficients $a_t$ and some factors $g_{k}^{x_k}$ with elliptic variable $x_k$. Since every such $g_{k}$ is elliptic with respect to~$X$, the following number is finite:
$$
C_k=\underset{z\in \mathbb{Z}}{\max} |g_k^z|_X.
$$
Therefore, for any values of elliptic variables (hidden in $u_1,\dots,u_s$), we have
$$
\overset{s}{\underset{j=1}{\sum}}|u_j|_X\leqslant \overset{n}{\underset{t=1}{\sum}}|a_t|_X+
\underset{k\in \{1,\dots,n\}\setminus \mathcal{M}}{\sum} C_k.
$$

\medskip

\noindent
We rewrite equation~\eqref{eq9} in the form
\begin{equation}\label{eq9a}
u_1'(g_{i_1}')^{x_{i_1}}\dots u_s'(g_{i_s}')^{x_{i_s}}=1,
\end{equation}
where $g_{i_j}'$ (and $h_{i_j}$) were defined in Step 2 and $u_j'=h_{i_{j-1}}^{-1}u_jh_{i_{j}}$, $j=1,\dots,s$, and addition is modulo $s$. We check that for any concrete values of elliptic variables, this equation satisfies assumptions of Proposition~\ref{bounded_coefficients}.

First we note that for any concrete values of elliptic variables, we have
\begin{equation}\label{eq9b}
\overset{s}{\underset{j=1}{\sum}}|u_j'|_X\leqslant \overset{n}{\underset{t=1}{\sum}}|a_t|_X+
\underset{k\in \{1,\dots,n\}\setminus \mathcal{M}}{\sum}C_k+2\overset{s}{\underset{j=1}{\sum}}|h_{i_j}|_X.
\end{equation}
Thus, the lengths of coefficients $u_j'$ in~\eqref{eq9a} are bounded from above by a constant, which is independent of values of elliptic variables. Second, by {\bf (d)}, every coefficient $g_{i_j}'$ in~\eqref{eq9a} belongs to one of the subgroups from the set $\{E_G(g_i)\,|\, i\in \mathcal{M}^{nc}\}$. Therefore, by Proposition~\ref{bounded_coefficients}, equation~\eqref{eq9a} is equivalent to a finite disjunction of finite systems of kind $\Phi_{(\sigma,b)}$, and each such system consists of
equations $R$ of the following two types:

($T_1$) $R$ is loxodromic with respect to $X$ and is over $E_G(g_{i})$ for some $i\in \mathcal{M}^{nc}$.

($T_2$) $R$ has the form $u'_{j_1}b_{\ell_1}^{-1}u'_{j_2}b_{\ell_2}^{-1}\dots u'_{j_t}b_{\ell_t}^{-1}=1$.
Since $u_1',\dots,u_s'$ contain only elliptic variables, $R$ is either elliptic with respect to $X$, or without variables.

\medskip

Suppose that $\Phi_{(\sigma,b)}$ contains equations without variables. If one of them, say $a=1$, where $a\in G$, is not valid in $G$, we replace $\Phi_{(\sigma,b)}$ by the equivalent system consisting of the single elliptic equation $a\cdot1^{x_1}1^{x_2}\dots 1^{x_n}=1$.
If all of them are identities in $G$, we delete them from $\Phi_{(\sigma,b)}$.
Denote the resulting system again by~$\Phi_{(\sigma,b)}$.

\medskip

{\it Step 4.} (Checking loxodromicity and ellipticity for the action of $G$ on $S$).
As noticed above, equations in $\Phi_{(\sigma,b)}$ are either loxodromic or elliptic with respect to~$X$. It remains to check that the same alternative holds with respect to the action $G\curvearrowright S$.

\medskip

{\bf Claim.}
For any $i\in \{1,\dots,n\}$, if $g_i\in {\rm Ell}(G\curvearrowright \Gamma)$, then $g_i\in {\rm Ell}(G\curvearrowright S)$.

\medskip

{\it Proof.}
An  equivalent formulation is: For any $i\in \{1,\dots,n\}$, if $g_i\in {\rm Lox}(G\curvearrowright S)$, then $g_i\in {\rm Lox}(G\curvearrowright \Gamma)$. In short: If $i\in \mathcal{L}$, then $g_i\in {\rm Lox}(G\curvearrowright \Gamma)$. The latter follows from $\mathcal{L}\subseteq \mathcal{M}$ and {\bf (b)}. \hfill $\Box$

\medskip


By this claim, equations in $\Phi_{(\sigma,b)}$ that are elliptic with respect to the action $G\curvearrowright \Gamma$ remain elliptic with respect to the action $G\curvearrowright S$.
Now suppose that $R$ is an equation in $\Phi_{(\sigma,b)}$ that is loxodromic with respect to $G\curvearrowright \Gamma$. By $(T_1)$-$(T_2)$,
the equation $R$ is over the virtually cyclic subgroup $E_G(g_i)$ for some $i\in \mathcal{M}^{nc}$.
If $g_i$ is loxodromic with respect to $G\curvearrowright S$, then $R$ is loxodromic with respect to $G\curvearrowright S$. If $g_i$ is elliptic with respect to $G\curvearrowright S$, then every element of $E_G(g_i)$ is elliptic with respect to $G\curvearrowright S$ as well. In this case $R$ is elliptic with respect to $G\curvearrowright S$.
\hfill $\Box$

\medskip

For the proof of Theorem~B, we need a version of Proposition~\ref{bounded_coefficients} without restriction on coefficients $a_i$.


\begin{proposition}\label{peripheral_equation}
Let $G$ be a group, $\{H_{\lambda}\}_{\lambda\in \Lambda}$ a collection of subgroups of~$G$, and $X$ a symmetric relative
generating set of $G$ with respect to $\{H_{\lambda}\}_{\lambda\in \Lambda}$ with $1\in X$.
Suppose that $\{H_{\lambda}\}_{\lambda\in \Lambda}$ is hyperbolically embedded in $G$ with respect to $X$.

Then any exponential equation of form \eqref{eq1}
$$
a_1g_1^{x_1}a_2g_2^{x_2}\dots a_ng_n^{x_n}=1
$$
\noindent
with $a_i\in G$ and $g_i\in H_{\lambda(i)}$ for some $\lambda(i) \in \Lambda$, $i=1,\dots,n$,
is equivalent to a finite disjunction $\Phi$ of finite conjunctions of exponential equations,
$$
\Phi:=\overset{k}{\underset{i=1}{\bigvee}}\overset{s_i}{\underset{j=1}{\bigwedge}} L_{i,j},
$$
such that equations in each conjunction are pairwise independent and each $L_{i,j}$ is an exponential equation over $H_{\lambda}$ for some $\lambda\in \Lambda$.

\end{proposition}

\medskip

{\it Proof.} We deduce this statement from Proposition~~\ref{bounded_coefficients} using the following trick. We write $a_i=a_{i,1}g_{i,1}\dots a_{i,t_i}g_{i,t_i}a_{i,t_i+1}$, where $a_{i,j}\in \langle X\rangle$ and $g_{i,j}\in H_{\lambda(i,j)}$ for some $\lambda(i,j)\in \Lambda$. Consider the exponential equation
$$
\overset{n}{\underset{i=1}{\prod}}a_{i,1}g_{i,1}^{x_{i,1}}\dots a_{i,t_i}g_{i,t_i}^{x_{i,t_i}}a_{i,{t_i+1}}g_i^{x_i}=1.
$$
By Proposition~\ref{bounded_coefficients} this equation is equivalent to
a finite disjunction $\Psi$ of finite conjunctions of exponential equations,
$$
\Psi:=\overset{k}{\underset{i=1}{\bigvee}}\overset{s_i}{\underset{j=1}{\bigwedge}} L_{i,j},
$$
such that equations in each conjunction are pairwise independent and each $L_{i,j}$ is an exponential equation over $H_{\lambda}$ for some $\lambda\in \Lambda$. Note that $\Psi$ depends on variables $x_{i,1},\dots x_{i,t_i},x_i$, where $i=1,\dots,n$. Then the desired $\Phi$ can be obtained from $\Psi$ by setting $x_{i,1}=\dots =x_{i,t_i}=1$ for $i=1,\dots,n$. \hfill $\Box$

\medskip

{\it Proof of Theorem B.} By definition of relative hyperbolicity, there exists a finite subset $X\subseteq G$ such that $G$ is generated by
$X$ together with all $H_{\lambda}$, $\lambda\in \Lambda$. We may assume that $X$ is symmetric and $1\in X$.
We set $
Y=X\sqcup (\underset{\lambda\in \Lambda}{\sqcup} H_{\lambda})
$.
Then the Cayley graph $S=\Gamma(G,Y)$ is hyperbolic and $G$ acts on $S$ acylindrically
(see~\cite[Corollary 2.54]{Osin_0} and~\cite[Proposition 5.2]{Osin_1}, respectively).
By Theorem~A, it suffices to consider the case where equation (1) is elliptic,
i.e., all $g_1,\dots,g_n$ are elliptic with respect to~$Y$.
Then every $g_i$ has finite order or is conjugate to an  element of one of the subgroups
$H_{\lambda}$, $\lambda\in \Lambda$ (see subsection~2.5).
Introducing finite number of finitary equations and using conjugations, we can reduce to the case where every $g_i$ lies in $H_{\lambda(i)}$ for some $\lambda(i)\in \Lambda$.
Finally, we apply Proposition~\ref{peripheral_equation}.
\hfill $\Box$

\section{Proofs of Corollaries}

\medskip

Statement (a) of Corollary C follows from Theorem A and the following two lemmas, where the first lemma is obvious.
Statement (b) of Corollary C follows from statement (a) and Theorem B.

\begin{lemma} If $A\subseteq \mathbb{Z}^s$ and $B\subseteq \mathbb{Z}^t$ are two $\mathbb{Z}$-semilinear subsets, 
then their concatenation $A\oplus B$ is a $\mathbb{Z}$-semilinear subset in $\mathbb{Z}^{s+t}$.
\end{lemma}

\begin{lemma} Let $G$ be a virtually cyclic group.
The solution set of an exponential equation \eqref{eq1}
with coefficients $a_1,\dots,a_n, g_1,\dots ,g_n, $ from $G$ is a $\mathbb{Z}$-semilinear subset of $\mathbb{Z}^n$.
\end{lemma}

\medskip

{\it Proof.} Let $H$ be some finite index normal cyclic subgroup of $G$. If $G$ is finite, we assume that $H=1$.
Let $N$ be the index of $H$ in $G$.
We write $x_i=Ny_i+z_i$, where $y_i$ and $z_i$ are new variables taking values in $\mathbb{Z}$ and $\mathcal{N}=\{0,\dots,N-1\}$, respectively. Then it suffices to prove that for each concrete tuple
$(k_1,\dots ,k_n)\in \mathcal{N}^n$ the solution set of the exponential equations
\begin{equation}\label{eq5}
a_1(g_1^N)^{y_1}g_1^{k_1}\dots a_n(g_n^N)^{y_n}g_n^{k_n}=1
\end{equation}
with variables $y_1,\dots,y_n$ is either empty or a $\mathbb{Z}$-linear subset of $\mathbb{Z}^n$.
If $H=1$, then $g_i^N=1$ for all $i$. In this case the solution set of~\eqref{eq5} is $\mathbb{Z}^n$ if $a_1g_1^{k_1}\dots a_ng_n^{k_n}=1$ and is $\emptyset$ otherwise. Therefore we assume that $H\neq 1$. Then $H$ is an infinite cyclic subgroup of $G$. Let $h$ be a generator of $H$.

Since $g_i^N\in H$, there is $s_i\in \mathbb{Z}$ such that $g_i^N=h^{s_i}$, $i=1,\dots,n$.
Using notation $f_i=a_1g_1^{k_1}\dots a_{i-1}g_{i-1}^{k_{i-1}}a_{i}$, we rewrite~\eqref{eq5} as
\begin{equation}\label{eq6}
\overset{n}{\underset{i=1}{\prod}} f_i h^{s_iy_i} f_i^{-1}=g_n^{-k_n}f_n^{-1}.
\end{equation}
Since $H$ is normal in $G$, there exist $\sigma_i\in \{-1,1\}$ such that $f_ihf_i^{-1}=h^{\sigma_i}$, and we rewrite \eqref{eq6} as
\begin{equation}\label{eq7}
\overset{n}{\underset{i=1}{\prod}} h^{\sigma_is_iy_i}=g_n^{-k_n}f_n^{-1}.
\end{equation}
We assume that $f_ng_n^{k_n}$ lies in $H$, otherwise the solution set is empty.
Then there exists $s\in \mathbb{Z}$ such that $g_n^{-k_n}f_n^{-1}=h^s$, and~\eqref{eq7} is equivalent to the equation
$$
\overset{n}{\underset{i=1}{\sum}} \sigma_is_iy_i=s.
$$
The solution set of this equation is either empty or $\mathbb{Z}$-linear.\hfill $\Box$

\medskip

Corollaries D and E follow straightforwardly from Corollary C and the following easy claims:

{\it Claim 1.}  The solution set (respectively, the $\mathbb{N}$-solution set) of any finitary exponential equation $ag^x=1$
over any group $G$ is $\mathbb{Z}$-semilinear (respectively, semilinear).


{\it Claim 2.} The set of $\mathbb{Z}$-semilinear (respectively, semilinear = definable in Presburger arithmetic, or definable in the weak Presburger arithmetic) sets is closed under taking of concatenation $\oplus$.
For every natural $n$, the set of definable subsets of $\mathbb{N}^n$ (respectively, of $\mathbb{Z}^n$)
in Presburger arithmetic (respectively, in the weak Presburger arithmetic) is boolean closed.

\def\refname{REFERENCES}

\end{document}